\documentclass{amsart}

\newtheorem{theorem}{Theorem}[section]
\newtheorem{lemma}[theorem]{Lemma}

\theoremstyle{definition}
\newtheorem{definition}[theorem]{Definition}

\theoremstyle{remark}
\newtheorem{remark}[theorem]{Remark}

\numberwithin{equation}{section}

\begin{document}

\title{\text{Higher-order Hahn's quantum variational calculus}}

\author[{A. M. C. Brito da Cruz}]{Artur M. C. Brito da Cruz}

\address{Escola Superior de Tecnologia de Set\'{u}bal,
Estefanilha, 2910-761 Set\'{u}bal, Portugal}

\email{artur.cruz@estsetubal.ips.pt}

\thanks{This work is part of the first author's PhD,
which is carried out at the University of Aveiro under the
Doctoral Programme \emph{Mathematics and Applications}
of Universities of Aveiro and Minho. Submitted
30-Sep-2010; revised 4-Jan-2011; accepted 19-Jan-2011;
for publication in \emph{Nonlinear Analysis Series A: 
Theory, Methods \& Applications}.}


\author[N. Martins]{Nat\'{a}lia Martins}

\address{Department of Mathematics,
University of Aveiro, 3810-193 Aveiro, Portugal}

\email{natalia@ua.pt}


\author[D. F. M. Torres]{Delfim F. M. Torres}

\address{Department of Mathematics,
University of Aveiro, 3810-193 Aveiro, Portugal}

\email{delfim@ua.pt}


\subjclass[2000]{Primary 39A13; Secondary 49K05}

\keywords{$q$-differences, Hahn's calculus, Euler-Lagrange equations}


\begin{abstract}
We prove a necessary optimality condition of
Euler--Lagrange type for quantum variational problems
involving Hahn's derivatives of higher-order.
\end{abstract}

\maketitle


\section{Introduction}

Many physical phenomena are described
by equations involving nondifferentiable functions,
\textrm{e.g.}, generic trajectories of quantum mechanics \cite{Feynman}.
Several different approaches to deal with nondifferentiable functions
are followed in the literature of variational calculus, including
the time scale approach, which typically deal with delta or nabla differentiable
functions \cite{Rui,malina,5}, the fractional approach, allowing to consider
functions that have no first order derivative but have fractional derivatives
of all orders less than one \cite{MyID:182,El-Nabulsi,gastao},
and the quantum approach, which is particularly useful to model physical
and economical systems \cite{Bang04,Cresson,4}.

Roughly speaking, a quantum calculus substitute the classical derivative
by a difference operator, which allows to deal with sets of nondifferentiable functions.
Several dialects of quantum calculus are available \cite{Ernst,Kac}.
For motivation to study a nondifferentiable quantum variational calculus
we refer the reader to \cite{Ric:Holder,Bang04,Cresson}.

In 1949 Hahn introduced the difference operator $D_{q,\omega}$ defined by
\[
D_{q,\omega}\left[  f\right]  \left(  t\right)  :=\frac{f\left(  qt+\omega
\right)  -f\left(  t\right)  }{\left(  q-1\right)  t+\omega},
\]
where $f$ is a real function, and $q\in (0,1)$ and $\omega>0$ are real fixed numbers \cite{Hahn}.
The Hahn difference operator has been applied successfully in the construction
of families of ortogonal polynomials as well as in approximation problems
\cite{alvares,odzi,Petronilho}.
However, during 60 years, the construction of the proper inverse
of Hahn's difference operator remained an open question.
Eventually, the problem was solved in 2009 by Aldwoah \cite{1}
(see also \cite{Aldwoah:IJMS,2}). Here we introduce
the higher-order Hahn's quantum variational calculus,
proving the Hahn quantum analog of the higher-order Euler--Lagrange equation.
As particular cases we obtain the $q$-calculus Euler--Lagrange equation
\cite{Bang04} and the $h$-calculus Euler--Lagrange equation
\cite{MyID:179,book:DCV}.

Variational functionals that depend on higher
derivatives arise in a natural way
in applications of engineering, physics, and economics.
Let us consider, for example, the equilibrium of an elastic
bending beam. Let us denote by $y(x)$ the deflection of the point $x$ of the beam,
$E(x)$ the elastic stiffness of the material, that can vary with $x$,
and $\xi(x)$ the load that bends the beam.
One may assume that, due to some constraints of physical nature, the dynamics
does not depend on the usual derivative $y'(x)$ but on some quantum derivative
$D_{q,\omega}\left[y\right]\left(x\right)$. In this condition,
the equilibrium of the beam correspond to the solution of the following
higher-order Hahn's quantum variational problem:
\begin{equation}
\label{ex:intro:ar}
\int_0^L \left[\frac{1}{2} \left(E(x) D_{q,\omega}^{2}\left[y\right](x)\right)^2
- \xi(x) y\left(q^2 x + q \omega + \omega\right)\right] dx \longrightarrow \min.
\end{equation}
Note that we recover the classical problem of the equilibrium
of the elastic bending beam when $(\omega,q)\rightarrow (0,1)$.
Problem \eqref{ex:intro:ar} is a particular case of the
problem \eqref{P} investigated in Section~\ref{sec:mr}.
Our higher-order Hahn's quantum Euler--Lagrange equation
(Theorem~\ref{Higher order E-L}) gives the main tool to solve such problems.

The paper is organized as follows. In Section~\ref{sec:prelim}
we summarize all the necessary definitions and properties
of the Hahn difference operator and the associated
$q,\omega$-integral. In Section~\ref{sec:mr} we formulate and prove our main results:
in \S\ref{sub:sec:hoCV} we prove a
higher-order fundamental Lemma of the calculus of variations
with the Hahn operator (Lemma~\ref{ordem n});
in \S\ref{sub:sec:HOELHCV} we deduce a higher-order
Euler--Lagrange equation for Hahn's variational calculus
(Theorem~\ref{Higher order E-L}); finally we provide in \S\ref{subsec:Ex}
a simple example of a quantum optimization problem
where our Theorem~\ref{Higher order E-L} leads
to the global minimizer, which is not a continuous function.


\section{Preliminaries}
\label{sec:prelim}

Let $q\in\left(0,1\right)$ and $\omega>0$. We introduce the real number
\[
\omega_{0}:=\frac{\omega}{1-q}.
\]
Let $I$ be a real interval containing $\omega_{0}$.
For a function $f$ defined on $I$,
the \emph{Hahn difference operator} of $f$ is given by
\[
D_{q,\omega}\left[f\right]\left(t\right)
:=\frac{f\left(  qt+\omega\right) - f\left(t\right)}{\left(q-1\right)t + \omega}\, , \text{ if }
t\neq\omega_{0}\, ,
\]
and $D_{q,\omega}\left[f\right]\left(\omega_{0}\right) := f^{\prime}\left(\omega_{0}\right)$,
provided $f$ is differentiable at $\omega_{0}$. We sometimes call $D_{q,\omega}\left[f\right]$
the $q,\omega$-\emph{derivative of} $f$, and $f$ is said to be $q,\omega$-\emph{differentiable on} $I$
if $D_{q,\omega}\left[f\right]\left(\omega_{0}\right)$ exists.

\begin{remark}
The $D_{q,\omega}$ operator generalizes (in the limit) the forward
$h$-difference and the Jackson $q$-difference operators \cite{Ernst,Kac}.
Indeed, when $q\rightarrow 1$ we obtain the forward $h$-difference
\[
\Delta_{h}\left[f\right]\left(t\right)
:=\frac{f\left(t+h\right) - f\left(t\right)}{h},
\]
when $\omega\rightarrow 0$ we obtain the Jackson $q$-difference operator
\[
D_{q}\left[f\right]\left(t\right)
:= \frac{f\left(qt\right)-f\left(t\right)}{t\left(  q-1\right)}\, ,
\text{ if }t\neq 0 \, ,
\]
and $D_{q}[f]\left(0\right) =f^{\prime}\left(0\right)$ provided
$f^{\prime}\left(0\right)$ exists. Notice also that,
under appropriate conditions,
\[
\lim_{\omega\rightarrow0,q\rightarrow1}D_{q,\omega}\left[f\right]\left(t\right)
=f^{\prime}\left(t\right).
\]
\end{remark}
The Hahn difference operator has the following properties:

\begin{theorem}[\cite{1,Aldwoah:IJMS,2}]
Let $f$ and $g$ be $q,\omega$-differentiable on $I$ and $t\in I$. One has:
\begin{enumerate}
\item  $D_{q,\omega}[f](t) \equiv 0$ on $I$ if and only if $f$ is constant;

\item $D_{q,\omega}\left[  f+g\right]  \left(  t\right)  =D_{q,\omega}\left[
f\right]  \left(  t\right)  +D_{q,\omega}\left[  g\right]  \left(  t\right)$;

\item $D_{q,\omega}\left[  fg\right]  \left(  t\right)  =D_{q,\omega}\left[
f\right]  \left(  t\right)  g\left(  t\right)  +f\left(  qt+\omega\right)
D_{q,\omega}\left[  g\right]  \left(  t\right)$;

\item $\displaystyle D_{q,\omega}\left[  \frac{f}{g}\right]  \left(  t\right)
=\frac{D_{q,\omega}\left[  f\right]  \left(  t\right)  g\left(  t\right)
-f\left(  t\right)  D_{q,\omega}\left[  g\right]  \left(  t\right)  }{g\left(
t\right)  g\left(  qt+\omega\right)  }$ if $g\left(
t\right)  g\left(  qt+\omega\right)  \neq0$;

\item $f\left(  qt+\omega\right)  =f\left(  t\right)  +\left(  t\left(
q-1\right)  +\omega\right)  D_{q,\omega}\left[  f\right]  \left(  t\right)  $.
\end{enumerate}
\end{theorem}

For $k\in\mathbb{N}_{0}=\mathbb{N}\cup\left\{  0\right\}$
define $\displaystyle\left[  k\right]  _{q}:=\frac{1-q^{k}}{1-q}$
and let $\sigma\left(t\right)  =qt+\omega$, $t\in I$.
Note that $\sigma$ is a contraction, $\sigma(I)\subseteq I$,
$\sigma\left(t\right)<t$ for $t>\omega_{0}$,
$\sigma\left(  t\right)>t$ for $t<\omega_{0}$,
and $\sigma\left(  \omega_{0}\right)  =\omega_{0}$.
The following technical result is used several times in our paper:

\begin{lemma}[\cite{1,2}]
Let $k\in\mathbb{N}$ and $t\in I$. Then,
\begin{enumerate}
\item
$\sigma^{k}\left(  t\right)  =\underset{k\text{-times}}{\underbrace{\sigma
\circ\sigma\circ \cdots \circ\sigma}}\left(  t\right)  =q^{k}t+\omega\left[
k\right]_{q}$;

\item
$\displaystyle \left(  \sigma^{k}\left(  t\right)  \right)^{-1}
=\sigma^{-k}\left(t\right)=\frac{t-\omega\left[  k\right]  _{q}}{q^{k}}$.
\end{enumerate}
\end{lemma}

>From now on $I$ denotes an interval of $\mathbb{R}$ containing $\omega_{0}$.
Following \cite{1,Aldwoah:IJMS,2} we define the notion of
$q,\omega$-\emph{integral} (also known as the \emph{Jackson--N\"{o}rlund integral})
as follows:

\begin{definition}
Let  $a,b\in I$ and $a<b$. For $f:I\rightarrow\mathbb{R}$
the $q,\omega$-\emph{integral of} $f$ from $a$ to $b$ is given by
\[
\int_{a}^{b}f\left(  t\right)  d_{q,\omega}t:=\int_{\omega_{0}}^{b}f\left(
t\right)  d_{q,\omega}t-\int_{\omega_{0}}^{a}f\left(  t\right)  d_{q,\omega
}t\text{,}
\]
where
\[
\int_{\omega_{0}}^{x}f\left(  t\right)  d_{q,\omega}t:=\left(  x\left(
1-q\right)  -\omega\right)  \sum_{k=0}^{+\infty}q^{k}f\left(  xq^{k}
+\omega\left[  k\right]  _{q}\right)  \text{, }x\in I\, ,
\]
provided that the series converges at $x=a$ and $x=b$. In that case, $f$ is
called $q,\omega$-\emph{integrable on} $\left[a,b\right]$. We say that
$f$ is $q,\omega$-\emph{integrable over} $I$ if it is $q,\omega$-\emph{integrable}
over $[a,b]$ for all $a,b\in I$.
\end{definition}

\begin{remark}
The $q,\omega$-\emph{integral} generalizes (in the limit) the Jackson
$q$-integral and the N\"{o}rlund's sum \cite{Kac}.
When $\omega\rightarrow0$, we obtain
the Jackson $q$-integral
\[
\int_{a}^{b}f\left(  t\right)  d_{q}t:=\int_{0}^{b}f\left(  t\right)
d_{q}t-\int_{0}^{a}f\left(  t\right)  d_{q}t\text{,}
\]
where
\[
\int_{0}^{x}f\left(  t\right)  d_{q}t:=x\left(  1-q\right)  \sum_{k=0}
^{+\infty}q^{k}f\left(  xq^{k}\right)  \text{.}
\]
When $q\rightarrow1$, we obtain the N\"{o}rlund's sum
\[
\int_{a}^{b}f\left(  t\right)  \Delta_{\omega}t :=\int_{+\infty}^{b}f\left(
t\right)  \Delta_{\omega}t-\int_{+\infty}^{a}f\left(  t\right)  \Delta
_{\omega}t,
\]
where
\[
\int_{+\infty}^{x}f\left(  t\right)  \Delta_{\omega}t :=-\omega\sum
_{k=0}^{+\infty}f\left(  x+k\omega\right)  \text{.}
\]
\end{remark}

It can be shown  that if $f:I\rightarrow\mathbb{R}$
is continuous at $\omega_{0}$,
then $f$ is $q,\omega$-\emph{integrable over} $I$ \cite{1,Aldwoah:IJMS,2}.

\begin{theorem}[Fundamental Theorem of Hahn's Calculus \cite{1,2}]
Assume that $f:I\rightarrow \mathbb{R}$ is continuous at $\omega_{0}$ and,
for each $x\in I$, define
\[
F\left(  x\right)  :=\int_{\omega_{0}}^{x}f\left(  t\right)  d_{q,\omega
}t\text{.}
\]
Then $F$ is continuous at $\omega_{0}$. Furthermore, $D_{q,\omega}\left[
F\right]  \left(  x\right)$ exists for every $x\in I$ with
$D_{q,\omega}\left[  F\right]  \left(  x\right)  =f\left(  x\right)$.
Conversely,
$\int_{a}^{b}D_{q,\omega}\left[  f\right]  \left(  t\right)  d_{q,\omega
}t=f\left(  b\right) -f\left(  a\right)$  for all $a,b\in I$.
\end{theorem}

The $q,\omega$-\emph{integral} has the following properties:

\begin{theorem}[\cite{1,Aldwoah:IJMS,2}]
\label{Propriedades do integral}
Let $f,g:I\rightarrow\mathbb{R}$ be $q,\omega$-\emph{integrable on} $I$,
$a,b,c\in I$ and $k\in\mathbb{R}$. Then,
\begin{enumerate}
\item $\int_{a}^{a}f\left(  t\right)  d_{q,\omega}t=0$;

\item $\int_{a}^{b}kf\left(  t\right)  d_{q,\omega}t=k\int_{a}^{b}f\left(
t\right)  d_{q,\omega}t$;

\item \label{eq:item3} $\int_{a}^{b}f\left(  t\right)  d_{q,\omega}t=-\int_{b}^{a}f\left(
t\right)  d_{q,\omega}t$;

\item $\int_{a}^{b}f\left(  t\right)  d_{q,\omega}t=\int_{a}^{c}f\left(
t\right)  d_{q,\omega}t+\int_{c}^{b}f\left(  t\right)  d_{q,\omega}t$;

\item $\int_{a}^{b}\left(  f\left(  t\right)  +g\left(  t\right)  \right)
d_{q,\omega}t=\int_{a}^{b}f\left(  t\right)  d_{q,\omega}t+\int_{a}
^{b}g\left(  t\right)  d_{q,\omega}t$;

\item Every Riemann integrable function $f$ on $I$
is $q,\omega$-\emph{integrable on} $I$;

\item \label{itm:ip} If $f,g:I\rightarrow \mathbb{R}$
are $q,\omega$-differentiable and $a,b\in I$, then
\[
\int_{a}^{b}f\left(  t\right)  D_{q,\omega}\left[  g\right]  \left(  t\right)
d_{q,\omega}t= f\left(  t\right)  g\left(  t\right)  \bigg|_{a}^{b}
-\int_{a}^{b}D_{q,\omega}\left[  f\right]  \left(  t\right)  g\left(
qt+\omega\right)  d_{q,\omega}t.
\]
\end{enumerate}
\end{theorem}

Property~\ref{itm:ip} of Theorem~\ref{Propriedades do integral} is known as
$q,\omega$-\emph{integration by parts}. Note that
\[
\int_{\sigma\left(  t\right)  }^{t}f\left(  \tau\right)  d_{q,\omega}
\tau=\left(  t\left(  1-q\right)  -\omega\right)  f\left(  t\right)  \text{.}
\]

\begin{lemma}[\textrm{cf.} \cite{1,2}]
\label{positividade}
Let $b \in I$ and $f$ be $q,\omega$-\emph{integrable} over $I$.
Suppose that
$$
f(t)\geq 0 \quad \forall
t\in\left\{  q^{n}b+\omega\left[  n\right]  _{q}:n\in \mathbb{N}_{0}\right\}.
$$
\begin{enumerate}

\item If $\omega_0 \leq b$, then
$$
\int_{\omega_0}^b f(t)d_{q,\omega}t\geq 0.
$$

\item If $\omega_0 > b$, then
$$
\int_b^{\omega_0} f(t)d_{q,\omega}t\geq 0.
$$
\end{enumerate}
\end{lemma}

\begin{remark}
There is an inconsistency in \cite{1,2}.
Indeed, Lemma~6.2.7 of \cite{1} is only valid
if $b \ge \omega_0$ and $a \le b$.
Similarly with respect
to Lemma~3.7 of \cite{2}.
\end{remark}

\begin{remark}
\label{rem:diff:int}
In general it is not true that
$$
\left\vert \int_{a}^{b}f\left(  t\right)  d_{q,\omega}t\right\vert
\leq \int_{a}^{b}| f\left(  t\right)|  d_{q,\omega}t , \ \ \ a,b \in I.
$$
For a counterexample see \cite{1,2}. This illustrates well the difference
with other non-quantum integrals, \textrm{e.g.}, the time scale integrals
\cite{MR2562284,australia}.
\end{remark}

For $s\in I$ we define
\begin{equation}
\label{eq:def:tar}
\left[  s\right]_{q,\omega}:=\left\{  q^{n}s+\omega\left[  n\right]_q  :n\in
\mathbb{N}_{0}\right\}  \cup\left\{  \omega_{0}\right\}  \text{.}
\end{equation}

The following definition and lemma are important
for our purposes.

\begin{definition}
Let $s \in I$ and $g:I\times(-\bar{\theta},\bar{\theta}) \rightarrow \mathbb{R}$.
We say that $g\left(t,\cdot\right)$ is differentiable at $\theta_{0}$ uniformly
in $\left[s\right]_{q,\omega}$ if for every $\varepsilon>0$ there exists
$\delta>0$ such that
\[
0<\left\vert \theta-\theta_{0}\right\vert <\delta
\Rightarrow
\left\vert \frac{g\left(  t,\theta\right)  -g\left(  t,\theta_{0}\right)
}{\theta-\theta_{0}}-\partial_{2} g\left(  t,\theta_{0}\right)  \right\vert
<\varepsilon
\]
for all $t\in\left[  s\right]_{q,\omega}$,
where $\displaystyle\partial_{2}g=\frac{\partial g}{\partial\theta}$.
\end{definition}

\begin{lemma}[\textrm{cf.} \cite{4}]
\label{derivada do integral}
Let $s \in I$. Assume that
$g:I\times(-\bar{\theta},\bar{\theta}) \rightarrow \mathbb{R}$
is differentiable at $\theta_{0}$ uniformly in $\left[s\right]_{q,\omega}$,
and $\displaystyle\int_{\omega_{0}}^{s}\partial_{2}g\left(t,\theta_{0}\right)
d_{q,\omega}t$ exist. Then,
$$
G\left(  \theta\right):=\int_{\omega_{0}}^{s}g\left(t,\theta\right) d_{q,\omega}t,
$$
for $\theta$ near $\theta_{0}$, is differentiable at
$\theta_{0}$ with $G^{\prime}\left(  \theta_{0}\right)
=\displaystyle\int_{\omega_{0}}^{s}\partial_{2}g\left(t,\theta_{0}\right) d_{q,\omega}t$.
\end{lemma}


\section{Main Results}
\label{sec:mr}

We define the $q,\omega$-derivatives of higher-order in the usual way:
the $r$th $q,\omega$-derivative ($r \in \mathbb{N}$)
of $f:I\rightarrow \mathbb{R}$ is the function
$D_{q,\omega}^{r}[f]: I\rightarrow \mathbb{R}$ given by
$D_{q,\omega}^{r}[f]:=D_{q,\omega}[D_{q,\omega}^{r-1}[f]]$,
provided $D_{q,\omega}^{r-1}[f]$ is $q,\omega$-differentiable on $I$
and where $D_{q,\omega}^{0}[f]:=f$.

Let $a,b \in I$ and $a<b$. We introduce the linear space
$\mathcal{Y}^{r} = \mathcal{Y}^{r}\left(\left[a,b\right],\mathbb{R}\right)$ by
$$
\mathcal{Y}^{r} :=
\left\{  y: I  \rightarrow \mathbb{R}\, |\,
D_{q,\omega}^{i}[y], i = 0,\ldots, r,
\text{ are bounded on $[a,b]$ and continuous at } \omega_{0}\right\}
$$
endowed with the norm
$\left\Vert y\right\Vert_{r,\infty}:=\sum_{i=0}^{r}\left\Vert D_{q,\omega}
^{i}\left[  y\right]  \right\Vert_{\infty}$,
where
$\left\Vert y\right\Vert_{\infty}:=\sup_{t\in\left[  a,b\right]  }\left\vert
y\left(  t\right)  \right\vert$. The following notations are in order:
$\sigma(t)=qt+\omega$, $y^{\sigma}(t) = y^{\sigma^1}(t)
= (y \circ \sigma)(t) = y\left(  qt+\omega\right)$,
and $y^{\sigma^{k}} = y \circ y^{\sigma^{k-1}}$, $k = 2, 3, \ldots$
Our main goal is to establish necessary optimality conditions
for the higher-order $q,\omega$-variational
problem\footnote{In problem \eqref{P} ``extr'' denotes
``extremize'' (\textrm{i.e.}, minimize or maximize).}
\begin{equation}
\label{P} \tag{P}
\begin{gathered}
\mathcal{L}\left[  y\right] =  \int_{a}^{b}L\left(
t,y^{\sigma^{r}}\left(  t\right)  ,D_{q,\omega}\left[  y^{\sigma^{r-1}}\right]
\left(  t\right)  ,\ldots,D_{q,\omega}^{r}\left[  y\right]  \left(  t\right)
\right)  d_{q,\omega}t \longrightarrow \textrm{extr}\\
y \in \mathcal{Y}^{r}\left(  \left[  a,b\right]  ,\mathbb{R}\right)\\
y\left(  a\right)  =\alpha_{0} \, , \quad  y\left(  b\right)  =\beta_{0} \, ,\\
\vdots\\
D_{q,\omega}^{r-1}\left[  y\right]  \left(  a\right)  =\alpha_{r-1} \, , \quad
D_{q,\omega}^{r-1}\left[  y\right]  \left(  b\right)  =\beta_{r-1}\, ,
\end{gathered}
\end{equation}
where $r\in \mathbb{N}$ and $\alpha_{i}, \beta_{i}\in
\mathbb{R}$, $i=0,\ldots,r-1$, are given.
\begin{definition}
We say that $y$ is an admissible function for \eqref{P}
if $y \in\mathcal{Y}^{r}\left(  \left[  a,b\right], \mathbb{R}\right)$
and $y$ satisfies the boundary conditions
$D_{q,\omega}^{i}\left[  y\right]  \left(  a\right)  =\alpha_{i}$
and $D_{q,\omega}^{i}\left[  y\right]  \left(  b\right)  =\beta_{i}$
of problem \eqref{P}, $i = 0,\ldots,r-1$.
\end{definition}

The Lagrangian $L$ is assumed to satisfy the following hypotheses:
\begin{enumerate}
\item[(H1)] $(u_0, \ldots, u_r)\rightarrow L(t,u_0, \ldots, u_r)$
is a $C^1(\mathbb{R}^{r+1}, \mathbb{R})$ function for any $t \in [a,b]$;

\item[(H2)] $t \rightarrow L(t, y(t), D_{q,\omega}\left[  y\right](t),
\ldots, D_{q,\omega}^{r}\left[  y\right](t))$ is continuous
at $\omega_0$ for any admissible $y$;

\item[(H3)]  functions $t \rightarrow \partial_{i+2}L(t, y(t),
D_{q,\omega}\left[  y\right](t), \cdots, D^{r}_{q,\omega}\left[  y\right](t))$,
$i=0, 1, \cdots, r$, belong to
$\mathcal{Y}^{1}\left(  \left[  a,b\right]  ,\mathbb{R}\right)$
for all admissible $y$.
\end{enumerate}

\begin{definition}
We say that $y_{\ast}$ is a local minimizer (resp. local maximizer) for problem
\eqref{P} if $y_{\ast}$ is an admissible function and there exists $\delta>0$ such that
\[
\mathcal{L}\left[  y_{\ast}\right]  \leq\mathcal{L}\left[  y\right]
\text{ \ \ (resp. }\mathcal{L}\left[  y_{\ast}\right]  \geq
\mathcal{L}\left[  y\right]  \text{) }
\]
for all admissible $y$ with
$\left\Vert y_{\ast}-y\right\Vert_{r,\infty}<\delta$.
\end{definition}

\begin{definition}
We say that $\eta\in \mathcal{Y}^{r}\left(  \left[  a,b\right],
\mathbb{R}\right) $ is a \emph{variation}
if $\eta\left(  a\right)  =\eta\left(  b\right)  =0$, \ldots,
$D_{q,\omega}^{r-1}\left[  \eta\right]  \left(  a\right)
=D_{q,\omega}^{r-1}\left[  \eta\right]  \left(  b\right)  =0$.
\end{definition}

We define the $q,\omega$-interval from $a$ to $b$ by
\[
\left[  a,b\right]_{q,\omega}:=\left\{  q^{n}a+\omega\left[  n\right]
_{q}:n\in
\mathbb{N}_{0}\right\}  \cup
\left\{  q^{n}b+\omega\left[  n\right]_{q}:n\in\mathbb{N}_{0}\right\}
\cup\left\{  \omega_{0}\right\},
\]
\textrm{i.e.}, $\left[a,b\right]_{q,\omega}=[a]_{q,\omega}\cup [b]_{q,\omega}$,
where $[a]_{q,\omega}$ and $[b]_{q,\omega}$ are given by \eqref{eq:def:tar}.


\subsection{Higher-order fundamental lemma of Hahn's variational calculus}
\label{sub:sec:hoCV}

The chain rule, as known from classical calculus,
does not hold in Hahn's quantum context
(see a counterexample in \cite{1,2}).
However, we can prove the following.

\begin{lemma}
\label{q}
If $f$ is $q,\omega$-differentiable on $I$, then the following equality holds:
\[
D_{q,\omega}\left[  f^{\sigma}\right]  \left(  t\right)  =q\left(
D_{q,\omega}\left[  f\right]  \right)  ^{\sigma}\left(  t\right)  \text{,
\ }t\in I\text{.}
\]
\end{lemma}

\begin{proof}
For $t\neq\omega_{0}$ we have
\begin{equation*}
\left(  D_{q,\omega}\left[  f\right]  \right)  ^{\sigma}\left(  t\right)
=\frac{f\left(  q\left(  qt+\omega\right)  +\omega\right)  -f\left(
qt+\omega\right)  }{\left(  q-1\right)  \left(  qt+\omega\right)  +\omega}
  =\frac{f\left(  q\left(  qt+\omega\right)  +\omega\right)  -f\left(
qt+\omega\right)  }{q\left(  \left(  q-1\right)  t+\omega\right)  }
\end{equation*}
and
\begin{equation*}
D_{q,\omega}\left[  f^{\sigma}\right]  \left(  t\right)
=\frac{f^{\sigma}\left(  qt+\omega\right)
-f^{\sigma}\left(  t\right)  }{\left(  q-1\right) t+\omega}
=\frac{f\left(  q\left(  qt+\omega\right)  +\omega\right)
-f\left(qt+\omega\right)}{\left(  q-1\right)  t+\omega}.
\end{equation*}
Therefore,
$D_{q,\omega}\left[  f^{\sigma}\right]  \left(  t\right)  =q\left(
D_{q,\omega}\left[  f\right]  \right)  ^{\sigma}\left(  t\right)$.
If $t=\omega_{0}$, then $\sigma\left(  \omega_{0}\right)  =\omega_{0}$. Thus,
\begin{equation*}
\left(  D_{q,\omega}\left[  f\right]  \right)  ^{\sigma}\left(  \omega
_{0}\right) =\left(  D_{q,\omega}\left[  f\right]  \right)  \left(
\sigma\left(  \omega_{0}\right)  \right)
=\left(  D_{q,\omega}\left[  f\right]  \right)  \left(  \omega_{0}\right)
=f^{\prime}\left(  \omega_{0}\right)
\end{equation*}
and
$D_{q,\omega}\left[  f^{\sigma}\right]  \left(  \omega_{0}\right) =\left[
f^{\sigma}\right]  ^{\prime}\left(  \omega_{0}\right)
=f^{\prime}\left(  \sigma\left(
\omega_{0}\right)  \right)  \sigma^{\prime}\left(  \omega_{0}\right)
=qf^{^{\prime}}\left(  \omega_{0}\right)$.
\end{proof}

\begin{lemma}
\label{anula}
If $\eta\in \mathcal{Y}^{r}\left(  \left[  a,b\right]  ,\mathbb{R}\right)$ is such that
$D_{q,\omega}^{i}\left[  \eta\right]  \left(
a\right)  =0$ (resp. $D_{q,\omega}^{i}\left[  \eta\right]  \left(  b\right)
=0$) for all $i\in\left\{  0,1,\ldots,r\right\}  ,$ then $D_{q,\omega}
^{i-1}\left[  \eta^{\sigma}\right]  \left(  a\right)  =0$ (resp. $D_{q,\omega
}^{i-1}\left[  \eta^{\sigma}\right]  \left(  b\right)  =0$) for all $i\in\left\{
1,\ldots,r\right\}  $.
\end{lemma}

\begin{proof}
If $a=\omega_{0}$ the result is trivial (because $\sigma\left(  \omega
_{0}\right)  =\omega_{0}$). Suppose now that $a\neq\omega_{0}$ and fix
$i\in\left\{  1,\ldots,r\right\}  $. Note that
\[
D_{q,\omega}^{i}\left[  \eta\right]  \left(  a\right)  =\frac{\left(
D_{q,\omega}^{i-1}\left[  \eta\right]  \right)  ^{\sigma}\left(  a\right)
-D_{q,\omega}^{i-1}\left[  \eta\right]  \left(  a\right)  }{\left(
q-1\right)  a+\omega}\text{.}
\]
Since, by hypothesis, $D_{q,\omega}^{i}\left[  \eta\right]  \left(  a\right)
=0$ and $D_{q,\omega}^{i-1}\left[  \eta\right]  \left(  a\right)  =0$, then
$\left(  D_{q,\omega}^{i-1}\left[  \eta\right]  \right)  ^{\sigma}\left(
a\right)  =0$. Lemma~\ref{q} shows that
\[
\left(  D_{q,\omega}^{i-1}\left[  \eta\right]  \right)  ^{\sigma}\left(
a\right)  =\left(  \frac{1}{q}\right)  ^{i-1}D_{q,\omega}^{i-1}\left[
\eta^{\sigma}\right]  \left(  a\right).
\]
We conclude that
$D_{q,\omega}^{i-1}\left[  \eta^{\sigma}\right]  \left(  a\right)  =0$.
The case $t=b$ is proved in the same way.
\end{proof}

\begin{lemma}
\label{Lema F CV}
Suppose that $f \in \mathcal{Y}^{1}\left(\left[a,b\right],\mathbb{R}\right)$.
One has
\[
\int_{a}^{b}f\left(  t\right)  D_{q,\omega}\left[  \eta\right]  \left(
t\right)  d_{q,\omega}t=0
\]
for all functions $\eta\in \mathcal{Y}^{1}\left(  \left[  a,b\right]  ,\mathbb{R}\right)$ such that
$\eta\left(  a\right)  =\eta\left(  b\right)  =0$ if and only if
$f\left(  t\right)  =c$, $c\in \mathbb{R}$, for all $t\in\left[  a,b\right]_{q,\omega}$.
\end{lemma}

\begin{proof}
The implication ``$\Leftarrow$'' is obvious.
We prove ``$\Rightarrow$''. We begin noting that
\[
\underset{=0}{\underbrace{\int_{a}^{b}f\left(  t\right)  D_{q,\omega}\left[
\eta\right]  \left(  t\right)  d_{q,\omega}t}}=\underset{=0}{\underbrace
{f\left(  t\right)  \eta\left(  t\right)  \bigg|_{a}^{b}}}-\int_{a}
^{b}D_{q,\omega}\left[  f\right]  \left(  t\right)  \eta^{\sigma}\left( t
\right)  d_{q,\omega}t.
\]
Hence,
\[
\int_{a}^{b}D_{q,\omega}\left[  f\right]  \left(  t\right)  \eta\left(
qt+\omega\right)  d_{q,\omega}t=0
\]
for any $\eta\in \mathcal{Y}^{1}\left(  \left[  a,b\right]  ,
\mathbb{R}
\right)  $  such that $\eta\left(  a\right)  =\eta\left(  b\right)  =0$.
We need to prove that, for some $c\in\mathbb{R}$,
$f\left(  t\right)  =c$ for all $t\in\left[  a,b\right]  _{q,\omega}$,
that is,
$D_{q,\omega}\left[  f\right]  \left(  t\right)  =0$ for all $t\in\left[
a,b\right]  _{q,\omega}$.
Suppose, by contradiction, that there exists $p\in\left[  a,b\right]_{q,\omega}$
such that $D_{q,\omega}\left[  f\right]  \left(  p\right)  \neq 0$.

\noindent (1) If $p\neq\omega_{0}$, then
$p=q^{k}a+\omega\left[  k\right]_{q}$ or $p=q^{k}b+\omega\left[k\right]_{q}$
for some $k\in \mathbb{N}_{0}$. Observe that
$a\left(  1-q\right) -\omega$ and $b\left(  1-q\right)-\omega$
cannot vanish simultaneously.

\noindent (a) Suppose that $a\left(  1-q\right)  -\omega\neq0$ and $b\left(
1-q\right)  -\omega\neq0$. In this case we can assume, without loss of
generality, that
$p=q^{k}a+\omega\left[  k\right]_{q}$
and we can define
\begin{equation*}
\eta\left(  t\right) =
\begin{cases}
D_{q,\omega}\left[  f\right]  \left(  q^{k}a+\omega\left[  k\right]_{q}\right)
&  \text{ if } t=q^{k+1}a+\omega\left[  k+1\right]  _{q}\\
0 &  \text{ otherwise.}
\end{cases}
\end{equation*}
Then,
\begin{multline*}
\int_{a}^{b}D_{q,\omega}\left[  f\right]  \left(  t\right)  \cdot\eta\left(
qt+\omega\right)  d_{q,\omega}t\\
=-\left(  a\left(  1-q\right)
-\omega\right)  q^{k}D_{q,\omega}\left[  f\right]  \left(  q^{k}
a+\omega\left[  k\right]  _{q}\right)\cdot
D_{q,\omega}\left[  f\right]  \left(  q^{k}a+\omega\left[  k\right]_{q}\right) \neq 0,
\end{multline*}
which is a contradiction.

\noindent (b) If $a\left(  1-q\right)  -\omega\neq0$ and $b\left(  1-q\right)
-\omega=0$, then $b=\omega_{0}$. Since
$q^{k}\omega_{0}+\omega\left[  k\right]_{q}=\omega_{0}$ for all
$k\in \mathbb{N}_{0}$, then
$p\neq q^{k}b+\omega\left[  k\right]  _{q}$ $\forall k\in
\mathbb{N}_{0}$ and, therefore,
\[
p=q^{k}a+\omega\left[  k\right]_{q,\omega}\text{ for some }
k\in \mathbb{N}_{0}\text{.}
\]
Repeating the proof of $\left(  a\right)$ we obtain again a contradiction.

\noindent (c) If $a\left(  1-q\right)  -\omega=0$ and $b\left(  1-q\right)
-\omega\neq0$ then the proof is similar to $\left(  b\right)$.\\

\noindent (2) If $p=\omega_{0}$ then, without loss of generality, we can assume
$D_{q,\omega}\left[  f\right]  \left(  \omega_{0}\right)  >0$. Since
\[
\lim_{n\rightarrow+\infty}\left(  q^{n}a+\omega\left[  k\right]  _{q}\right)
=\lim_{n\rightarrow+\infty}\left(  q^{n}b+\omega\left[  k\right]  _{q}\right)
=\omega_{0}
\]
(see \cite{1}) and $D_{q,\omega}\left[  f\right]  $ is continuous at
$\omega_{0}$, then
\begin{equation*}
\lim_{n\rightarrow+\infty}D_{q,\omega}\left[  f\right]  \left(  q^{n}
a+\omega\left[  k\right]  _{q}\right)   =\lim_{n\rightarrow
+\infty}D_{q,\omega}\left[  f\right]
\left(  q^{n}b+\omega\left[  k\right]_{q}\right)
=D_{q,\omega}\left[  f\right]  \left(  \omega_{0}\right)>0.
\end{equation*}
Thus, there exists $N\in \mathbb{N}$ such that for all $n\geq N$
one has
$D_{q,\omega}\left[  f\right]  \left(  q^{n}a+\omega\left[  k\right]_{q}\right)  >0$
and $D_{q,\omega}\left[  f\right]  \left(  q^{n}b+\omega\left[  k\right]_{q}\right)  >0$.

\noindent (a) If $\omega_{0}\neq a$ and $\omega_{0}\neq b$, then we can define
\[
\eta\left(  t\right)  =\left\{
\begin{tabular}[c]{lll}
$D_{q,\omega}\left[  f\right]  \left(  q^{N}b+\omega\left[  N\right]
_{q}\right)  $ & if & $t=q^{N+1}a+\omega\left[  N+1\right]  _{q}$ \\
&  & \\
$D_{q,\omega}\left[  f\right]  \left(  q^{N}a+\omega\left[  N\right]
_{q}\right)  $ & if & $t=q^{N+1}b+\omega\left[  N+1\right]  _{q}$ \\
&  & \\
0 &  & \text{otherwise}.
\end{tabular}
\right.
\]
Hence,
\begin{multline*}
\int_{a}^{b}D_{q,\omega}\left[  f\right]  \left(  t\right)  \eta\left(
qt+\omega\right)  d_{q,\omega}t\\
=\left(  b-a\right)  \left(  1-q\right)q^{N}D_{q,\omega}\left[  f\right]
\left(  q^{N} b+\omega\left[ N\right]  _{q}\right)\cdot  D_{q\omega}\left[  f\right]
\left(  q^{N}a+\omega\left[  N\right]  _{q}\right) \neq 0,
\end{multline*}
which is a contradiction.

\noindent (b) If $\omega_{0}=b$, then we define
\[
\eta\left(  t\right)  =\left\{
\begin{tabular}
[c]{lll}
$D_{q,\omega}\left[  f\right]  \left(  \omega_{0}\right)  $ & if &
$t=q^{N+1}a+\omega\left[  N+1\right]  _{q}$\\
&  & \\
$0$ &  & otherwise.
\end{tabular}
\right.
\]
Therefore,
\begin{equation*}
\begin{split}
\int_{a}^{b} & D_{q,\omega}\left[  f\right]  \left(  t\right)  \eta\left(
qt+\omega\right)  d_{q,\omega}t\\
 &=  -\int_{\omega_{0}}^{a}D_{q,\omega}\left[
f\right]  \left(  t\right)  \eta\left(  qt+\omega\right)  d_{q,\omega}t\\
 &= -\left(  a\left(  1-q\right)  -\omega\right)   q^{N}D_{q,\omega
}\left[  f\right]  \left(  q^{N}a+\omega\left[  k\right]  _{q}\right)
\cdot D_{q,\omega}\left[  f\right]  \left(  \omega_{0}\right) \neq  0,
\end{split}
\end{equation*}
which is a contradiction.

\noindent (c) When $\omega_{0}=a$, the proof is similar to $\left(  b\right)$.
\end{proof}

\begin{lemma}[Fundamental lemma of Hahn's variational calculus]
\label{Ordem 1}
Let $f,g \in \mathcal{Y}^{1}\left(  \left[  a,b\right]  ,
\mathbb{R}\right).$
If
\[
\int_{a}^{b}\left(  f\left(  t\right)  \eta^{\sigma}\left(  t\right)
+g\left(  t\right)  D_{q,\omega}\left[  \eta\right]  \left(  t\right)
\right)  d_{q,\omega}t=0
\]
for all $\eta\in \mathcal{Y}^{1}\left(  \left[  a,b\right]  ,
\mathbb{R}\right) $  such that $\eta\left(  a\right)  =\eta\left(  b\right)  =0$, then
\[
D_{q,\omega}\left[  g\right]  \left(  t\right)  =f\left(  t\right)  \text{
\ }\forall t\in\left[  a,b\right]  _{q,\omega}\text{.}
\]
\end{lemma}

\begin{proof}
Define the function $A$  by $A\left(  t\right)
:=\int_{\omega_{0}}^{t}f\left(  \tau\right)  d_{q,\omega}\tau$. Then
$D_{q,\omega}\left[  A\right]  \left(  t\right)  =f\left(  t\right)  $ for all
$t\in\left[  a,b\right]  $\textbf{\ } and
\begin{align*}
\int_{a}^{b}A\left(  t\right)  D_{q,\omega}\left[  \eta\right]  \left(
t\right)  d_{q,\omega}t  & =  A\left(  t\right)  \eta\left(  t\right)
\bigg|_{a}^{b}-\int_{a}^{b}D_{q,\omega}\left[  A\right]  \left(  t\right)
\eta^{\sigma}\left(  t\right)  d_{q,\omega}t\\
& =-\int_{a}^{b}D_{q,\omega}\left[  A\right]  \left(  t\right)  \eta^{\sigma
}\left(  t\right)  d_{q,\omega}t\\
& =-\int_{a}^{b}f\left(  t\right)  \eta^{\sigma}\left(  t\right)  d_{q,\omega
}t\text{.}
\end{align*}
Hence,
\begin{multline*}
\int_{a}^{b}\left(  f\left(  t\right)  \eta^{\sigma}\left(  t\right)
+g\left(  t\right)  D_{q,\omega}\left[  \eta\right]  \left(  t\right)
\right)  d_{q,\omega}t=0\\
 \Leftrightarrow  \int_{a}^{b}\left(  -A\left(  t\right)  +g\left(  t\right)
\right)  D_{q,\omega}\left[  \eta\right]  \left(  t\right)  d_{q,\omega}t=0\text{.}
\end{multline*}
By Lemma~\ref{Lema F CV} there is a
$c\in \mathbb{R}$ such that $-A\left(  t\right)  +g\left(  t\right)  =c$ for all
$t\in\left[a,b\right]_{q,\omega}$. Hence $D_{q,\omega}\left[  A\right]  \left(
t\right)  =D_{q,\omega}\left[  g\right]  \left(  t\right)  $ for $t\in\left[
a,b\right]_{q,\omega}$, which provides the desired result:
$D_{q,\omega}\left[  g\right]  \left(  t\right)  =f\left(  t\right)  \text{
\ }\forall t\in\left[  a,b\right]_{q,\omega}$.
\end{proof}

We are now in conditions to deduce the
higher-order fundamental Lemma of
Hahn's quantum variational calculus.

\begin{lemma}[Higher-order fundamental lemma of Hahn's variational calculus]
\label{ordem n}
Let $f_{0},f_{1},\ldots,f_{r} \in \mathcal{Y}^{1}\left(  \left[  a,b\right]  ,\mathbb{R}\right)$. If
\[
\int_{a}^{b}\left(  \sum_{i=0}^{r}f_{i}\left(  t\right)  D_{q,\omega}
^{i}\left[  \eta^{\sigma^{r-i}}\right]  \left(  t\right)  \right)
d_{q,\omega}t=0
\]
for any variation $\eta$, then
\[
\sum_{i=0}^{r}\left(  -1\right)  ^{i}\left(  \frac{1}{q}\right)
^{\frac{\left(  i-1\right)  i}{2}}D_{q,\omega}^{i}\left[  f_{i}\right]
\left(  t\right)  =0
\]
for all $t\in\left[  a,b\right]_{q,\omega}$.
\end{lemma}

\begin{proof}
We proceed by mathematical induction. If $r=1$ the result is true
by Lemma~\ref{Ordem 1}. Assume that
\[
\int_{a}^{b}\left(  \sum_{i=0}^{r+1}f_{i}\left(  t\right)  D_{q,\omega}^{i}
\left[  \eta^{\sigma^{r+1-i}}\right]  \left(  t\right)  \right)  d_{q,\omega
}t=0
\]
for all functions $\eta$ such that $\eta\left(  a\right)
= \eta\left(b\right)  =0$, \ldots,
$D_{q,\omega}^{r}\left[  \eta\right]  \left(  a\right)
= D_{q,\omega}^{r}\left[  \eta\right]  \left(  b\right)  =0$. Note that
\begin{align*}
\int_{a}^{b} f_{r+1}&\left(  t\right)  D_{q,\omega}^{r+1}\left[  \eta\right]
\left(  t\right)  d_{q,\omega}t\\
 &=f_{r+1}\left(  t\right)
D_{q,\omega}^{r}\left[  \eta\right]  \left(  t\right)  \bigg|_{a}^{b}
 -\int_{a}^{b}D_{q,\omega}\left[  f_{r+1}\right]  \left(  t\right)  \left(
D_{q,\omega}^{r}\left[  \eta\right]  \right)  ^{\sigma}\left(  t\right)
d_{q,\omega}t\\
 &= -\int_{a}^{b}D_{q,\omega}\left[  f_{r+1}\right]  \left(  t\right)  \left(
D_{q,\omega}^{r}\left[  \eta\right]  \right)  ^{\sigma}
\left(  t\right)  d_{q,\omega}t
\end{align*}
and, by Lemma~\ref{q},
\[
\int_{a}^{b}f_{r+1}\left(  t\right)  D_{q,\omega}^{r+1}\left[  \eta\right]
\left(  t\right)  d_{q,\omega}t=-\int_{a}^{b}D_{q,\omega}\left[
f_{r+1}\right]  \left(  t\right)  \left(  \frac{1}{q}\right)  ^{r}D_{q,\omega
}^{r}\left[  \eta^{\sigma}\right]  \left(  t\right)  d_{q,\omega}t\text{.}
\]
Therefore,
\begin{align*}
\int_{a}^{b}&\left(  \sum_{i=0}^{r+1}f_{i}\left(  t\right)  D_{q,\omega}^{i}
\left[  \eta^{\sigma^{r+1-i}}\right]  \left(  t\right)  \right)  d_{q,\omega
}t \\
&=\int_{a}^{b}\left(  \sum_{i=0}^{r}f_{i}\left(  t\right)  D_{q,\omega
}^{i}\left[  \eta^{\sigma^{r+1-i}}\right]  \left(  t\right)
\right)  d_{q,\omega}t\\
&\qquad -\int_{a}^{b}D_{q,\omega}\left[  f_{r+1}\right]  \left(  t\right)  \left(
\frac{1}{q}\right)  ^{r}D_{q,\omega}^{r}\left[  \eta^{\sigma}\right]  \left(
t\right)  d_{q,\omega}t\\
&=\int_{a}^{b}
\biggl[
\sum_{i=0}^{r-1}f_{i}\left(  t\right)  D_{q,\omega}^{i}\left[  \left(
\eta^{\sigma}\right)  ^{\sigma^{r-i}}\right]  \left(  t\right)  d_{q,\omega
}t\\
& \qquad +\left(  f_{r}-\left(  \frac{1}{q}\right)  ^{r}D_{q,\omega}\left[
f_{r+1}\right]  \right)  \left(  t\right)  D_{q,\omega}^{r}\left[
\eta^{\sigma}\right]  \left(  t\right)
\biggr]
d_{q,\omega}t.
\end{align*}
By Lemma~\ref{anula}, $\eta^{\sigma}$ is a variation.
Hence, using the induction hypothesis,
\begin{align*}
\sum_{i=0}^{r-1}&\left(-1\right)^{i}\left(  \frac{1}{q}\right)^{\frac{\left(
i-1\right)  i}{2}}D_{q,\omega}^{i}\left[  f_{i}\right]
\left(  t\right)\\
& \qquad +\left(  -1\right)  ^{r}\left(  \frac{1}{q}\right)^{\frac{\left(
r-1\right)  r}{2}}D_{q,\omega}^{r}\left[  \left(  f_{r}
-\frac{1}{q^{r}}D_{q,\omega}\left[  f_{r+1}\right]  \right)  \right]  \left(
t\right)\\
&=\sum_{i=0}^{r-1}\left(  -1\right)  ^{i}\left(  \frac{1}{q}\right)^{\frac{\left(
i-1\right)  i}{2}}D_{q,\omega}^{i}\left[  f_{i}\right]
\left(  t\right)  +\left(  -1\right)  ^{r}\left(  \frac{1}{q}\right)^{\frac{\left(
r-1\right)  r}{2}}D_{q,\omega}^{r}\left[  f_{r}\right]
\left(  t\right) \\
& \qquad +\left(  -1\right)  ^{r+1}\left(  \frac{1}{q}\right)  ^{\frac{\left(
r-1\right)  r}{2}}\frac{1}{q^{r}}D_{q,\omega}^{r}\left[  D_{q,\omega}\left[
f_{r+1}\right]  \right]  \left(  t\right) \\
&=0
\end{align*}
for all $t\in\left[  a,b\right]_{q,\omega}$, which leads to
\[
\sum_{i=0}^{r+1}\left(  -1\right)  ^{i}\left(  \frac{1}{q}\right)
^{\frac{\left(  i-1\right)  i}{2}}D_{q,\omega}^{i}\left[  f_{i}\right]
\left(  t\right)  =0\text{, \ }t\in\left[  a,b\right]  _{q,\omega}.
\]
\end{proof}


\subsection{Higher-order Hahn's quantum Euler--Lagrange equation}
\label{sub:sec:HOELHCV}

For a variation $\eta$ and an admissible function $y$, we define the function
$\phi:\left(-\bar{\epsilon},\bar{\epsilon}\right)  \rightarrow \mathbb{R}$ by
\[
\phi\left(  \epsilon\right)  =\phi\left(  \epsilon,y,\eta\right)
:=\mathcal{L}\left[  y+\epsilon\eta\right]  .
\]
The first variation of the variational problem \eqref{P} is defined by
\[
\delta\mathcal{L}\left[  y,\eta\right]  :=
\phi^{\prime}\left(  0\right)  .
\]
Observe that
\begin{align*}
\mathcal{L}\left[  y+\epsilon\eta\right]
=&\int_{a}^{b}L\Bigg(t,y^{\sigma^{r}}\left(  t\right)
+\epsilon\eta^{\sigma^{r}}\left(  t\right)
,D_{q,\omega}\left[  y^{\sigma^{r-1}}\right]  \left(  t\right)  +\epsilon
D_{q,\omega}\left[  \eta^{\sigma^{r-1}}\right]  \left(  t\right),\\
& \qquad\qquad\qquad \ldots,
D_{q,\omega}^{r}\left[  y\right]  \left(  t\right)
+\epsilon D_{q,\omega}^{r}\left[
\eta\right]  \left(  t\right)  \Bigg)d_{q,\omega}t\\
&=\mathcal{L}_{b}\left[y+\epsilon\eta\right]
-\mathcal{L}_{a}\left[  y+\epsilon\eta\right]
\end{align*}
with
\begin{multline*}
\mathcal{L}_{\xi}\left[  y+\epsilon\eta\right]
=\int_{\omega_{0}}^{\xi}L\Bigg(t,y^{\sigma^{r}}\left(  t\right)
+\epsilon\eta^{\sigma^{r}}\left(t\right), D_{q,\omega}\left[
y^{\sigma^{r-1}}\right]  \left(  t\right)
+\epsilon D_{q,\omega}\left[  \eta^{\sigma^{r-1}}\right]  \left(  t\right),\\
\ldots, D_{q,\omega}^{r}\left[  y\right]  \left(  t\right)
+\epsilon D_{q,\omega}^{r}\left[  \eta\right]  \left(  t\right)  \Bigg)d_{q,\omega}t,
\end{multline*}
$\xi \in \{a, b\}$. Therefore,
\begin{equation}
\label{eq:le}
\delta\mathcal{L}\left[  y,\eta\right]  =\delta\mathcal{L}_{b}\left[
y,\eta\right]  -\delta\mathcal{L}_{a}\left[  y,\eta\right]  .
\end{equation}
Considering \eqref{eq:le}, the following lemma is a direct
consequence of Lemma~\ref{derivada do integral}:

\begin{lemma}
\label{uniformly}
For a variation $\eta$ and an admissible function $y$, let
\begin{multline*}
g\left(  t,\epsilon\right) :=L\bigg(t,y^{\sigma^{r}}\left(  t\right)
+\epsilon\eta^{\sigma^{r}}\left(  t\right)  ,D_{q,\omega}\left[
y^{\sigma^{r-1}}\right]  \left(  t\right)  +\epsilon D_{q,\omega}\left[
\eta^{\sigma^{r-1}}\right]  \left(  t\right),\\
\ldots, D_{q,\omega}^{r}\left[  y\right]  \left(  t\right)
+\epsilon D_{q,\omega}^{r}\left[  \eta\right]  \left(  t\right)  \bigg),
\end{multline*}
$\epsilon\in\left(-\bar{\epsilon},\bar{\epsilon}\right)$. Assume that:

\noindent (1) $g\left(  t,\cdot\right)$ is differentiable at $0$ uniformly in
$t\in \left[a,b\right]_{q,\omega}$;

\noindent (2) $\mathcal{L}_{a}\left[  y+\epsilon\eta\right]
=\displaystyle\int_{\omega_{0}}^{a}g\left(
t,\epsilon\right)  d_{q,\omega}t $ and $\mathcal{L}_{b}\left[
y+\epsilon\eta\right]=\displaystyle\int_{\omega_{0}}^{b}g\left(
t,\epsilon\right)  d_{q,\omega}t$ exist for $\epsilon \approx 0$;

\noindent (3) $\displaystyle\int_{\omega_{0}}^{a}\partial_{2}g\left(  t,0\right)
d_{q,\omega}t$ and $\displaystyle\int_{\omega_{0}}^{b}\partial_{2}g\left(
t,0\right)  d_{q,\omega}t$ exist.

\noindent Then
\begin{multline*}
\phi^{\prime}\left(  0\right) = \delta\mathcal{L}\left[  y,\eta\right]
= \int_{a}^{b}\bigg(  \sum_{i=0}^{r}\partial_{i+2}L\left(t,
y^{\sigma^{r}}\left(  t\right), D_{q,\omega}\left[
y^{\sigma^{r-1}}\right]  \left(  t\right),
\ldots, D_{q,\omega}^{r}\left[  y\right]  \left(  t\right)\right)\\
\cdot D_{q,\omega}^{i}\left[
\eta^{\sigma^{r-i}}\right]  \left(  t\right)\bigg) d_{q,\omega}t,
\end{multline*}
where $\partial_{i}L$ denotes the partial derivative
of $L$ with respect to its $i$th argument.
\end{lemma}

The following result gives a necessary condition of Euler--Lagrange type
for an admissible function to be a local extremizer for \eqref{P}.

\begin{theorem}[Higher-order Hahn's quantum Euler--Lagrange equation]
\label{Higher order E-L}
Under hypotheses (H1)--(H3) and conditions (1)--(3)
of Lemma~\ref{uniformly} on the Lagrangian $L$, if
$y_{\ast}\in \mathcal{Y}^{r}$ is a local extremizer for problem \eqref{P},
then $y_{\ast}$ satisfies the $q,\omega$-Euler--Lagrange equation
\begin{equation}
\label{eq:E-L}
\sum_{i=0}^{r}\left(  -1\right)  ^{i}\left(  \frac{1}{q}\right)^{\frac{\left(
i-1\right)  i}{2}}D_{q,\omega}^{i}\left[  \partial_{i+2} L\right]\left(
t,y^{\sigma^{r}}\left(  t\right), D_{q,\omega}\left[
y^{\sigma^{r-1}}\right]  \left(  t\right),
\ldots, D_{q,\omega}^{r}\left[  y\right]\left(t\right)\right)
= 0
\end{equation}
for all $t\in\left[  a,b\right]_{q,\omega}$.
\end{theorem}

\begin{proof}
Let $y_{\ast}$ be a local extremizer for problem
\eqref{P} and $\eta$ a variation. Define
$\phi:(-\bar{\epsilon},\bar{\epsilon}) \rightarrow \mathbb{R}$
by $\phi\left(  \epsilon\right)  :=\mathcal{L}\left[  y_{\ast}+\epsilon\eta\right]$.
A necessary condition for $y_{\ast}$ to be an extremizer is given by
$\phi^{\prime}\left(  0\right)  =0$. By Lemma~\ref{uniformly} we conclude that
\begin{multline*}
\int_{a}^{b}\Biggl(  \sum_{i=0}^{r}\partial_{i+2}
L\left(  t,y^{\sigma^{r}}\left(  t\right)  ,D_{q,\omega}\left[
y^{\sigma^{r-1}}\right]  \left(  t\right),
\ldots, D_{q,\omega}^{r}\left[y\right] \left(  t\right)  \right)\\
\cdot D_{q,\omega}^{i}\left[\eta^{\sigma^{r-i}}\right]\left(
t\right)  \Biggr)  d_{q,\omega}t=0
\end{multline*}
and \eqref{eq:E-L} follows from Lemma~\ref{ordem n}.
\end{proof}

\begin{remark}
In practical terms the hypotheses of Theorem~\ref{Higher order E-L}
are not so easy to verify \emph{a priori}. One can, however,
assume that all hypotheses are satisfied and apply the
$q,\omega$-Euler--Lagrange equation \eqref{eq:E-L} heuristically
to obtain a \emph{candidate}. If such a candidate is, or not,
a solution to problem \eqref{P} is a different question that always requires
further analysis (see an example in \S\ref{subsec:Ex}).
\end{remark}

When $\omega\rightarrow 0$ one obtains from \eqref{eq:E-L}
the higher-order $q$-Euler--Lagrange equation:
\[
\sum_{i=0}^{r}\left(  -1\right)  ^{i}\left(  \frac{1}{q}\right)
^{\frac{\left(  i-1\right)  i}{2}}D_{q}^{i}\left[  \partial_{i+2}
L\right]  \left(  t,y^{\sigma^{r}}\left(  t\right)  ,D_{q}\left[
y^{\sigma^{r-1}}\right]  \left(  t\right),
\ldots, D_{q}^{r}\left[  y\right]
\left(  t\right)  \right)    =0
\]
for all $t\in\left\{  aq^{n}:n\in
\mathbb{N}_{0}\right\}  \cup\left\{  bq^{n}:n\in
\mathbb{N}_{0}\right\}  \cup\left\{  0\right\}$.
The higher-order $h$-Euler--Lagrange equation
is obtained from \eqref{eq:E-L} taking the limit $q\rightarrow 1$:
\[
\sum_{i=0}^{r}\left(  -1\right)  ^{i}
\Delta_{h}^{i}\left[  \partial_{i+2} L\right]\left(t,
y^{\sigma^{r}}\left(  t\right),\Delta_{h}\left[
y^{\sigma^{r-1}}\right]  \left(  t\right),
\ldots, \Delta_{h}^{r}\left[  y\right]  \left(  t\right)\right)=0
\]
for all $t\in\left\{  a+nh:n\in
\mathbb{N}_{0}\right\}
\cup\left\{  b+nh:n\in \mathbb{N}_{0}\right\}$.
The classical Euler--Lagrange equation \cite{6}
is recovered when $(\omega, q)\rightarrow (0, 1)$:
$$
\sum_{i=0}^{r}\left(-1\right)^{i}
\frac{d^i}{d t^i}\partial_{i+2}L\left(t,
y\left(t\right),y'\left(t\right), \ldots, y^{(r)}(t)\right)=0
$$
for all $t \in [a,b]$.

We now illustrate the usefulness
of our Theorem~\ref{Higher order E-L}
by means of an example that is not covered
by previous available results in the literature.


\subsection{An Example}
\label{subsec:Ex}

Let $q=\frac{1}{2}$ and $\omega=\frac{1}{2}$. Consider the following problem:
\begin{equation}
\label{example}
\mathcal{L}\left[  y\right]
=  \int_{-1}^{1}\left(y^\sigma(t) + \frac{1}{2}\right)^2  \left(
\left( D_{q,\omega}\left[  y\right] (t)\right)^2 -1 \right)^2 d_{q,\omega}t
\longrightarrow \min
\end{equation}
over all $y \in \mathcal{Y}^1$ satisfying the boundary conditions
\begin{equation}
\label{eq:bc:ex}
y(-1)=0 \ \ \ \text{and} \ \ \ y(1)=-1.
\end{equation}
This is an example of problem \eqref{P} with $r=1$.
Our $q,\omega$-Euler--Lagrange equation \eqref{eq:E-L}
takes the form
\[
D_{q,\omega}\left[\partial_{3}L\right]\left(
t,y^{\sigma}\left(  t\right)  ,D_{q,\omega}\left[  y\right]  \left(t\right)  \right)
=\partial_{2}L\left(  t,y^{\sigma} \left(  t\right)  ,D_{q,\omega}\left[
y\right]  \left(  t\right)  \right).
\]
Therefore, we look for an admissible function $y_{\ast}$
of \eqref{example}-\eqref{eq:bc:ex} satisfying
\begin{multline}
\label{eq:EL:ex}
D_{q,\omega}\left[4 \left(y^\sigma
+ \frac{1}{2}\right)^2\left( \left( D_{q,\omega}\left[
y\right]\right)^2 -1 \right)D_{q,\omega}\left[  y\right] \right](t)\\
= 2\left(y^\sigma(t) + \frac{1}{2}\right) \left(
\left( D_{q,\omega}\left[  y\right](t)\right)^2 -1 \right)
\end{multline}
for all $ t\in\left[  -1,1\right]_{q,\omega}$. It is easy to see that
\begin{equation*}
y_{\ast}(t) =
\begin{cases}
-t & \text{ if } t \in (-1,0)\cup (0,1]\\
0 & \text{ if } t=-1\\
1 & \text{ if } t=0
\end{cases}
\end{equation*}
is an admissible function for \eqref{example}-\eqref{eq:bc:ex} with
\begin{equation*}
D_{q,\omega}\left[  y_{\ast}\right] (t) =
\begin{cases}
-1 & \text{ if } t \in (-1,0)\cup (0,1]\\
1 & \text{ if } t=-1\\
-3 & \text{ if }  t=0,
\end{cases}
\end{equation*}
satisfying the $q,\omega$-Euler--Lagrange equation \eqref{eq:EL:ex}.
We now prove that the \emph{candidate} $y_{\ast}$ is indeed a minimizer
for \eqref{example}-\eqref{eq:bc:ex}.
Note that here $\omega_0=1$ and, by Lemma~\ref{positividade}
and item \eqref{eq:item3} of Theorem~\ref{Propriedades do integral},
\begin{equation}
\label{eq:inq:crz}
\mathcal{L}\left[  y\right] = \int_{-1}^{1}\left(y^\sigma(t) + \frac{1}{2}\right)^2\left(
\left( D_{q,\omega}\left[  y\right] (t)\right)^2 -1 \right)^2 d_{q,\omega}t   \geq 0
\end{equation}
for all admissible functions
$y \in \mathcal{Y}^1\left(\left[-1,1\right],\mathbb{R}\right)$.
Since $\mathcal{L}\left[  y_{\ast}\right]=0$,
we conclude that $y_{\ast}$ is a minimizer
for problem \eqref{example}-\eqref{eq:bc:ex}.

It is worth to mention that the minimizer $y_{\ast}$ of \eqref{example}-\eqref{eq:bc:ex}
is not continuous while the classical calculus of variations \cite{6},
the calculus of variations on time scales \cite{Rui,malina,5},
or the nondifferentiable scale variational calculus \cite{Ric:Holder,Ric,Cresson},
deal with functions which are necessarily continuous.
As an open question, we pose the problem of determining
conditions on the data of problem \eqref{P} assuring, \emph{a priori},
the minimizer to be regular.


\section*{Acknowledgments}

The first author is supported by the \emph{Portuguese Foundation
for Science and Technology} (FCT) through the PhD fellowship
SFRH/BD/33634/2009; the second and third authors by FCT through
the \emph{Center for Research and Development
in Mathematics and Applications} (CIDMA).
The authors are very grateful to the referee
for valuable remarks and comments.



\end{document}